\documentclass[12pt]{article}
\usepackage{amssymb}
\usepackage[dvipdfm]{graphicx}
\usepackage{amsfonts}
\usepackage{latexsym}
\usepackage{amsthm}
\usepackage{amsmath}

\usepackage{amssymb, amscd}

\usepackage{url}
\usepackage{gclc}

\newtheorem{problem}{Problem}%[section]
\newtheorem{theo}[problem]{Theorem}
\newtheorem{rem}[problem]{Remark}
\newtheorem{prob}[problem]{Problem}
\newtheorem{defin}[problem]{Definition}
\newtheorem{prop}[problem]{Proposition}

\newtheorem{exam}[problem]{Example}
\newtheorem{observ}[problem]{Observation}
\newtheorem{conj}[problem]{Conjecture}

\begin{document}

 \title{{Fulton-MacPherson compactification,\\  cyclohedra, and the polygonal pegs problem}}

\author{Sini\v sa T.\ Vre\' cica\\ {\small Faculty of Mathematics}\\[-2mm] {\small University of Belgrade}
\\[-2mm] {\small vrecica$@$matf.bg.ac.yu}
 \and Rade  T.\ \v Zivaljevi\' c\\ {\small Mathematical Institute}\\[-2mm] {\small SANU, Belgrade}\\[-2mm]
 {\small rade$@$mi.sanu.ac.yu} }
\date{\textit{\small To Anders Bj\" orner,\\ on the occasion of his
60th anniversary.}} \maketitle \vspace{-1cm}
\begin{abstract}
The cyclohedron  $W_n$, known also as the Bott-Taubes polytope,
arises both as the polyhedral realization of the poset of all
cyclic bracketings of the word $x_1x_2\ldots x_n$ and as an
essential part of the Fulton-MacPherson compactification of the
configuration space of $n$ distinct, labelled points on the circle
$S^1$. The ``polygonal pegs problem'' asks whether every simple,
closed curve in the plane or in the higher dimensional space
admits an inscribed polygon of a given shape. We develop a new
approach to the polygonal pegs problem based on the
Fulton-MacPherson (Axelrod-Singer, Kontsevich) compactification of
the configuration space of (cyclically) ordered $n$-element
subsets in $S^1$. Among the results obtained by this method are
proofs of Gr\" unbaum's conjecture about affine regular hexagons
inscribed in smooth Jordan curves and a new proof of the
conjecture of Hadwiger about inscribed parallelograms in smooth,
simple, closed curves in the $3$-space (originally established by
Makeev in \cite{Mak}).
\end{abstract}

\renewcommand{\thefootnote}{$\ast$} \footnotetext{Supported by Grants 144014
and 144026 of the Serbian Ministry of Science and Technology.}
\section{Introduction}

The classical ``square peg problem'', going back to Toeplitz
(1911) and Emch (1913), asks whether every Jordan curve in the
plane has four points forming a square. In the first published
account of the problem \cite{Em} the result was established for
the case of closed convex curves. Over the span of almost one
hundred years many interesting cases of the problem were resolved,
occasionally after initial partial refutations and subsequent
improvements over the original proofs. The reader is referred to
\cite{Gru}, \cite{KlWa}, and \cite{Pa} for a more complete
overview of the history of the problem and a brief discussion of
the leading contributions by Emch (1911, 1913), Shnirelman (1929,
1944), Jerrard (1961), Stromquist (1989), Griffiths (1991), and
others.

In spite of its simple and mathematically attractive formulation
the square peg problem has not been solved so far in full
generality. Other problems of similar nature were formulated in
the meantime and some of them have remained unsolved even in the
case of smooth curves, see \cite{Gru} and \cite{h71} for examples
and \cite{Gri}, \cite{Pa} for a broader outlook on the whole area.

In this paper we emphasize the role of {\em cyclohedra} $W_n$ in
the ``square peg problem'' and other related problems of discrete
geometry where polygons are inscribed in curves, surfaces etc. We
start in Section~\ref{sec:square-pegs} with a complete, reasonably
short and conceptually transparent (if not entirely elementary)
solution of the ``square peg problem'' in the case of smooth
curves. More importantly, this section serves as a model example
of how the ``cyclohedron approach'' (or the method of canonical
compactifications) can be applied to other problems of this
nature. In Section~\ref{sec:Grunbaum}, using this method, we prove
the Gr\"{u}nbaum's conjecture about inscribed affine regular
hexagons in smooth, simple closed curves in the plane. By the same
technique we prove a result in Section~\ref{sec:Hadwiger}
(established earlier in \cite{Mak} by different methods) which
confirms a conjecture of Hadwiger about inscribed parallelograms
in smooth, simple, closed curves in the $3$-space. These results
should not be seen as isolated examples. Rather they are an
indication of the potential of the method for applications to many
other classes of problems where the degeneration of point
configurations and the appearance of the associated
pseudo-solutions has been one of the main obstacles for applying
standard topological methods. In Section~\ref{sec:conclusion} we
briefly discuss possibilities for extending results from
Sections~\ref{sec:square-pegs}, \ref{sec:Grunbaum} and
\ref{sec:Hadwiger} to larger classes of curves and offer a broader
outlook to the method of canonical compactifications
(\textit{FMASK}-compactifications).

\section{Outline of the main idea}
\label{sec:outline}

Given a Jordan curve $\Gamma\subset \mathbb{R}^2$ and its
parametrization $f : S^1\rightarrow\Gamma$, the {\em configuration
space} of all (labelled) quadrangles inscribed in $\Gamma$ is
parameterized by the torus $T^4\cong (S^1)^4$. In order to
determine which of these quadrangles are squares one introduces an
associated {\em test map} $\Phi : T^4\rightarrow U$ where $U\cong
\mathbb{R}^4$ is the associated {\em test vector space}. The test
map is well chosen if $q\in T^4$ is associated to a square
inscribed in $\Gamma$ if and only if $\Phi(q)=0\in U$.

\begin{itemize} \item Recall that the configuration space, the test
map, and the test space are the basic ingredients of the well
known ``configuration space/test map scheme'' \cite{Z04}
($CS/TM$-scheme for short) for applying (equivariant) topological
methods in combinatorial geometry. This proof scheme has been
applied for decades before it was codified and named in
\cite{User1}, \cite{User2} and remains one of the main tools for
applying topological methods in geometric combinatorics.
\end{itemize}
One of the main difficulties with the application of the
$CS/TM$-scheme in the ``square peg problem'' and its relatives is
the appearance of pseudo-solutions, i.e.\ degenerate
configurations which pass the test $\Phi(q)=0$ but are not genuine
solutions. Indeed, the test map $\Phi$ often takes into account
only the mutual distances of elements of the configuration $q$, so
for example in the square peg problem $F$ does not distinguish
degenerate squares $q=(v,v,v,v)$ from actual squares.

A natural way to get around this difficulty is to remove from the
configuration space $T^4$ the diagonal $\Delta =\{q\mid
q=(v,v,v,v) \mbox{ {\rm for some} } v\in S^1\}$ or perhaps more
consistently the ``fat diagonal'' $\Delta_f :=
\{q=(q_1,q_2,q_3,q_4)\mid q_i=q_j \mbox{ {\rm for some} } i\neq
j\}$. The resulting truncated configuration spaces
$T^4\setminus\Delta$ and $F(S^1,4):= T^4\setminus\Delta_f$ are no
longer compact and this is often a source of other difficulties of
topological nature.

\medskip
Our main new idea is to ``blow up'' the degenerate configurations
in $\Delta$ ($\Delta_f$) and to modify accordingly (regularize)
the test map $\Phi$. This means that we replace the original
configuration space by the Fulton-MacPherson compactification
\cite{FM} of the truncated configurations space. Actually we use
its spherical version and a close relative due to Axelrod-Singer
and Kontsevich \cite{AS} \cite{Ko}, here referred to as the
canonical or \textit{FMASK}-compactification of the configuration
space.

In the context of the square peg problem, the configuration space
$F(S^1,4):=(S^1)^4\setminus\Delta_f$ is compactified to the
associated \textit{FMASK} compactification $F[S^1,4]$. The new
test map $\Psi : F[S^1,4]\rightarrow U$ is defined as the
extension of the map $\Phi' : F(S^1,4)\rightarrow U$ for a
suitable modification $\Phi'$ of the original test map which
essentially takes into account the rescaling of the degenerate
configurations.

Throughout the paper we mainly work with the subspace
$S^1(n)\subset F(S^1,n):=(S^1)^n\setminus\Delta_f$ of all
cyclically ordered $n$-tuples of points in $S^1$ and the
corresponding compactification $S^1[n]\subset F[S^1,n]$. As it was
shown in \cite{BT}, $S^1[n]\cong S^1\times W_n$, where $W_n$ is a
close relative of Stasheff polytope (associahedron) called
cyclohedron. This allows us to give a direct and elementary
exposition of $F[S^1,n]$ and $S^1[n]$ which is sufficient for all
our applications and which is fairly independent of the general
theory of \textit{FMASK} compactifications (see however \cite{Si}
for a more complete treatment of $F[M,n]$ and other related
compactifications).

\section{\textit{FMASK}-compactification of configuration spaces}
\label{sec:fmask}

Let us recall some elementary facts about the partially ordered
set $(\mathcal{C}(Y), \leqslant)$ of all compactifications of a
(locally compact, Hausdorff) space $Y$.

\begin{enumerate}\item[$\bullet_1$]
A compactification of $Y$ is a pair $cY=(X,c)$ where $c :
Y\rightarrow X$ is a homeomorphic embedding and $Cl_X(c(Y))=X$. By
definition  $c_1Y\leqslant c_2Y$ if there exists a continuous map
$f : c_2Y\rightarrow c_1Y$ such that $c_1 = f\circ c_2$. Two
compactifications $c_1Y$ and $c_2Y$ are considered equivalent (and
often identified) if both $c_1Y\leqslant c_2Y$ and $c_2Y\leqslant
c_1Y$, which turns ``$\leqslant$'' into an order relation on the
set $\mathcal{C}(Y)$ of (equivalence classes) of compactifications
of $Y$. A simple but important fact is that each non-empty set
$\mathcal{C}_0\subset \mathcal{C}(Y)$ has the least upper bound
with respect to the order ``$\leqslant$'', see e.g.\ \cite{Eng},
Theorem~3.5.9. Indeed, if $\mathcal{C}_0 = \{c_jY\}_{j\in J}$ then
the smallest compactification $\tau Y$ such that $c_jY\leqslant
\tau Y$ for all $j\in J$ can be described as the closure of the
image of the diagonal embedding $\Delta : Y\rightarrow \prod_{j\in
J}c_jY$.

\item[$\bullet_2$] A map  $g : Y\rightarrow Z$ from $Y$ to a
compact $Z$ is not necessarily extendable to a compactification
$cY$ prescribed in advance. However, there exists the smallest
compactification $\tau Y$ such that $cY\leqslant \tau Y$ and a map
$g' : \tau Y\rightarrow Z$ such that $g = g'\circ \tau$. Indeed,
it is easy to show that the closure $Cl_{cY\times Z}(\Gamma(g))$
of the graph $\Gamma(g)\subset Y\times Z\subset cY\times Z$, in
the compact space $cY\times Z$, has all the required properties.
More generally, given a family $\mathcal{F}=\{g_j\}_{j\in J}$ of
maps $g_j : Y\rightarrow Z_j$, from $Y$ to compact spaces $Z_j$,
there is the smallest compactification $\tau
Y=\tau_{\mathcal{F}}Y$ greater than $cY$ where all functions $g_j$
can be extended.
\end{enumerate}

\begin{exam}\label{exam:primer} {\rm  Suppose that
$Y=(S^1)^2\setminus\Delta$ is the space of all ordered pairs of
distinct points in $S^1$. Let $cY=T^2=(S^1)^2$ be its ``naive
compactification'' and let $g : Y\rightarrow S^1$ be the map
defined by $g(x,y):=(x-y)/\Vert x-y\Vert$. Then the
compactification $\tau Y$, described as the closure of the graph
$\Gamma(g)=\{(p,g(p))\mid p=(x,y)\in (S^1)^2\setminus\Delta
\}\subset (S^1)^2\times S^1$, is the ``oriented blow up'' of
$(S^1)^2$ along the diagonal $\Delta$. For the future reference
(Section~\ref{sec:Hadwiger}) we denote this compactification by
$\tilde{F}[S^1,2]$ and observe that it is homeomorphic to the
annulus $S^1\times [0,1]$.}

\end{exam}

Our main example of the construction of $\tau_{\mathcal{F}}Y$ is
the canonical compactification (or the \textit{FMASK}
compactification) $S^1[n]$ of the configurations space
$Y:=S^1(n)\subset (S^1)^n$ of all $n$-element subsets
$q=\{q_1\prec q_2\prec\ldots\prec q_n\prec q_1 \}\subset S^1$ of
cyclically ordered points in $S^1$.

Given consecutive indices $i-1,i,i+1$ (where $n+1:=1$), let
$\mathcal{F}=\{\theta_i\}_{i=1}^n$ be the collection of functions
$\theta_i : S^1(n)\rightarrow [0,1]$ defined by $\theta_i(q):=
\measuredangle(q_{i-1}q_i)/ \measuredangle(q_{i-1}q_{i+1})$. The
functions $\theta_i$ alow us to reconstruct $q\in S^1(n)$, up to a
rotation. However these functions cannot be extended to the
closure $cY:=Cl(S^1(n))$ of $S^1(n)$ in $(S^1)^n$, since for
example $\theta_i(q)$ is meaningless if $q_i=q_{i+1}$. For this
reason it is quite natural to pass to the compactification
$\tau_{\mathcal{F}}Y$ where all these functions are well defined.

\begin{defin}\label{def:FMASK}
The canonical or \textit{FMASK} compactification $S^1[n]$ of the
space $Y:= S^1(n)$ of all cyclically ordered $n$-element
configurations in $S^1$ is the compactification $\tau
Y=\tau_{\mathcal{F}}Y$ associated to $cY:=Cl(S^1(n))\subset
(S^1)^n$ and the family $\mathcal{F}=\{\theta_j\}_{j=1}^n$. More
explicitly, $\tau Y$ is the closure of the image of the embedding
$$
\Lambda: Y\hookrightarrow cY \times \prod_{j=1}^n [0,1]_{(j)}
$$
where $\Lambda = c\times \prod_j\theta_j$ is the associated
diagonal map. Similarly, by starting with the configuration space
$F(S^1,n)$ of all (not necessarily cyclically ordered) $n$-tuples
of distinct points in $S^1$, one obtains the associated
\textit{FMASK} compactification $F[S^1,n]$.
\end{defin}

\begin{enumerate}\item[$\bullet_3$]
The construction of the compactification $\tau_{\mathcal{F}}Y$
depends functorially on the family $\mathcal{F}$. This means that
$\tau_{\mathcal{F}}\leqslant \tau_{\mathcal{F'}}$ if all functions
from $\mathcal{F}$ are (informally speaking) expressible by
functions from $\mathcal{F'}$. In particular it is not difficult
to formulate a criterion when two compactifications
$\tau_{\mathcal{F}}$ and $\tau_{\mathcal{F'}}$ are equivalent.
This can be used to show the equivalence of $S^1[n]$ and
$F[S^1,n]$ with the more general constructions of $F[M,n]$
developed in \cite{AS}, \cite{Ko}, \cite{Si}, see also
Definition~\ref{def:FMASK-2}.
\end{enumerate}
For completeness and as an additional illustration of the main
construction described in $\bullet_2$ we finish this section with
the definition of the \textit{FMASK}-compactification
$F[\mathbb{R}^d,n]$. Note that the word ``compactification'' is
not quite appropriate here, however the construction of the
(partial) compactification  $\tau_{\mathcal{F}}Y$ is still
meaningful and natural from the geometric point of view.
\begin{defin}\label{def:FMASK-2}
The canonical or \textit{FMASK}-compactification
$F[\mathbb{R}^d,n]$ of the space $Y=
F(\mathbb{R}^d,n):=(\mathbb{R}^d)^n\setminus \Delta_f$ of all
collections of $n$, distinct, labelled points in $\mathbb{R}^d$ is
the (partial) compactification $\tau Y=\tau_{\mathcal{F}}Y$
associated to the naive ``compactification'' $cY:=
(\mathbb{R}^d)^n$ and the family
$\mathcal{F}=\{\alpha_{ij}\}_{1\leq i<\leq n}\cup
\{\beta_{ijk}\}_{i<j<k}$ where $\alpha_{ij}: Y\rightarrow S^{d-1}$
is defined by $\alpha_{ij}(q):=(q_j-q_i)/\parallel
q_j-q_i\parallel$ while $\beta_{ijk}: Y\rightarrow [0,+\infty]$ is
the function that records the ratio $\beta_{ijk}(q):=\parallel
q_i-q_j\parallel/\parallel q_i-q_k\parallel$.
\end{defin}

\section{Cyclohedron $W_n$}\label{sec:cyclohedron}

The following proposition reveals the stratified manifold
structure of the space $S^1[n]$.

\begin{theo}\label{thm:cyclo}{\rm (\cite{BT})} For $n\geq 3$,
\begin{equation}\label{eqn:cyclo}
S^1[n] \cong S^1\times W_n \end{equation} where $W_n$ is a
$(n-1)$-dimensional, convex polytope, called cyclohedron or the
Bott-Taubes polytope. $W_n$ is combinatorially described as the
convex polytope whose face lattice is isomorphic to the poset of
all partial, cyclic bracketings of the word $x_1x_2\ldots x_n$.
\end{theo}

\begin{figure}[hbt]
\centering
\includegraphics[scale=.50]{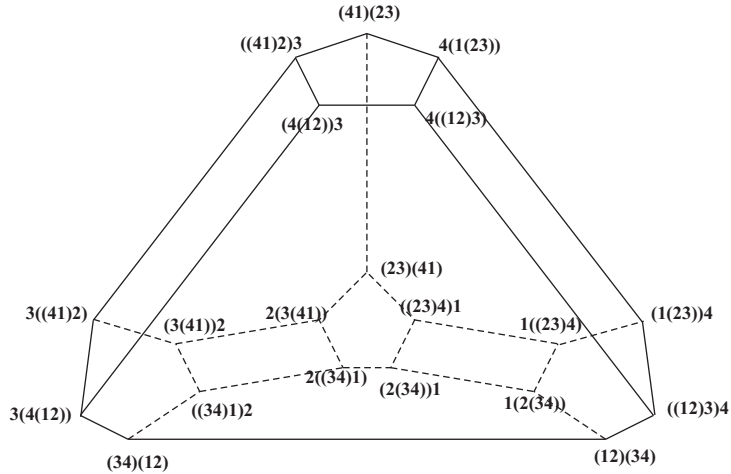}
\caption{Cyclohedron $W_4$.} \label{fig:cyclohedron}
\end{figure}

\medskip\noindent{\bf Proof:} (outline)
The reader is referred to \cite{BT}, \cite{Ma}, \cite{MSS}, and
\cite{Si} for more detailed presentation and related background
facts. We restrict ourselves to a brief explanation of the
isomorphism (\ref{eqn:cyclo}), sufficient for intended
applications.

The functions $\theta_i : S^1(n)\rightarrow [0,1]$ and their
extensions $\bar\theta_i : S^1[n]\rightarrow [0,1]$ can be used as
``coordinate functions'' on  spaces $S^1(n)$ and $S^1[n]$
respectively. They can be combined to create $2$-dimensional,
$3$-dimensional, or higher dimensional ``navigation instruments'',
with the corresponding screens being one, two, or higher
dimensional simplices $\Delta^2, \Delta^3$ etc. For example, given
a $4$-element subconfiguration $q_i\prec q_j\prec q_k\prec q_l$ of
$q=\{q_1\prec\ldots\prec q_n\}$, one can extend the function
$\lambda : S^1(n)\rightarrow \Delta^2$ defined by
\begin{equation}\label{eqn:2-screen}
\lambda(q):=\frac{1}{\measuredangle(q_i,q_l)}(\measuredangle(q_i,q_j),
\measuredangle(q_j,q_k),\measuredangle(q_k,q_l))
\end{equation}
to a function $\bar\lambda : S^1[n]\rightarrow \Delta^2$, where
$\measuredangle(p,q)=\widehat{p\,q}$ is the arc length of the
(counterclockwise) arc with endpoints $p$ and $q$. Indeed, the
function $\lambda$ can be expressed in terms of functions
$\theta_i$, consequently it can be extended to $S^1[n]$ and its
extension similarly expressed in terms of functions
$\bar\theta_i$.

\begin{figure}[hbt]
\centering
\includegraphics[scale=.50]{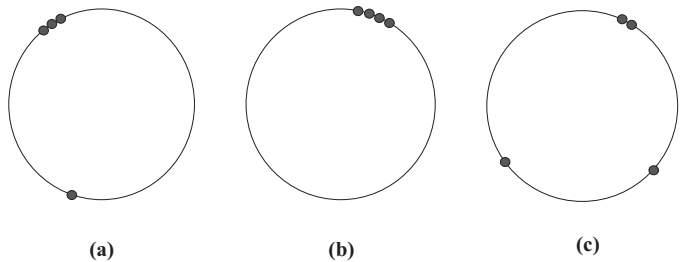}
\caption{Generic configurations in $\partial W_4$.}
\label{fig:3-circles}
\end{figure}

The configuration space $S^1(n)$ is clearly isomorphic to
$S^1\times {\rm Int}(\Delta)^{n-1}$. The reader can use the
``navigation screens'' to convince herself that the
compactification of this space is indeed described by equation
(\ref{eqn:cyclo}). For example one can check that the generic
configurations depicted in Figure~\ref{fig:3-circles} (a), (b),
and (c), respectively correspond to parallelograms, pentagonal,
and hexagonal facets of the cyclohedron $W_4$. \hfill $\square$

\begin{figure}[hbt]
\centering
\includegraphics[scale=.90]{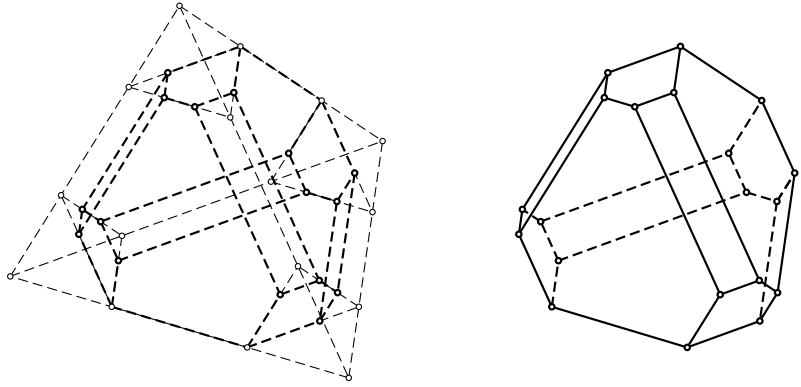}
\vspace{-2cm}\caption{Cyclohedron as a compactification/truncation
of a simplex.} \label{fig:cyclo}
\end{figure}

\section{Square pegs in round holes}
\label{sec:square-pegs}

We begin with a version of the ``square pegs in round holes''
theorem for $C^1$-smooth curves embedded in the $2$-space. This
result was in this generality proved by Stromquist \cite{St}, see
also Schnirelmann \cite{Shn} and Guggenheimer \cite{Gug} for
earlier results established with some extra hypotheses on the
smoothness or curvature of the curve. The reader is referred to
\cite{Gru} (p.\ 84) for a list of references addressing the case
of a convex curve and to \cite{Pak} for what appears to be the
only elementary presentation of the case of simple closed
polygons.

\begin{theo}\label{thm:smooth-square-pegs}
Every simple closed curve $\Gamma\subset \mathbb{R}^2$, which is
$C^1$-smooth, i.e.\ has a non-vanishing and continuously moving
tangent vector at each point, admits an inscribed square.
\end{theo}

Our proof of Theorem~\ref{thm:smooth-square-pegs} serves as a
model for other proofs of similar nature.  For this reason it is
broken into relevant individual steps illustrating
\textit{FMASK}-compactification modification of the usual
$CS/TM$-scheme \cite{Z04}. The scheme of the proof is summarized
in Section~\ref{subsec:proof-smooth-sqpegs}.

Theorem~\ref{thm:smooth-square-pegs} is a not the most general
result about inscribed squares in Jordan curves, see \cite{St},
and  Section~\ref{sec:conclusion} for a related discussion from
the \textit{FMASK}-compactification point of view.

\subsection{Configuration space $S^1(4)$ and the test maps}
\label{subsec:config}

Suppose that $f : S^1\rightarrow \mathbb{R}^2$ is a smooth
embedding satisfying the conditions of
Theorem~\ref{thm:smooth-square-pegs}. Moreover we silently assume,
here and elsewhere in the paper, that the embedding is
``counterclockwise'' in the sense that the degree of the
associated map $s\mapsto df/ds$ is $+1$.

Let $\Gamma = {\rm Image}(f)\subset \mathbb{R}^2$ be the
associated smooth curve. Suppose that  $U:= U_1\oplus U_2 \oplus
U_3$ is a vector space such that $U_1\cong \mathbb{R}^2, U_2\cong
U_3\cong\mathbb{R}$ and let $\Phi : (S^1)^4\rightarrow U$ be the
map defined by
\begin{equation}\label{eqn:peg-test-map}
\Phi(t) = (\phi_1(F(t)),\phi_2(F(t)),\phi_3(F(t)))
\end{equation}
where $F :=f^{\times 4} : (S^1)^4\rightarrow \mathbb{R}^2\oplus
\mathbb{R}\oplus \mathbb{R}$ is the map induced by $f$ and
$$
 \begin{array}{l}
\phi_1(y)=\frac{y_1+y_3}{2}-\frac{y_2+y_4}{2}, \\
\phi_2(y)=\|y_1-y_3 \| - \|y_2-y_4\|, \\
\phi_3(y)=\|y_1-y_2\| - \|y_2-y_3\| + \|y_3-y_4\| - \|y_4-y_1\|.
 \end{array}
$$
It is clear that $f(t_1),f(t_2),f(t_3),f(t_4)$ are consecutive
vertices of a square inscribed in the curve $\Gamma$ if and only
if $\Phi(t)=0$.

Let $\Phi_0 : S^1(4) \rightarrow U$ be the restriction of $\Phi$
on the configuration space $S^1(4)$ of all labelled $4$-element
subsets of $S^1$ such that the labelling agrees with the
counterclockwise (cyclic) order of points on the circle $S^1$.

The symmetric group $S_4$ acts on $(S^1)^4$ by permuting
coordinates. However, it is its subgroup $\mathbb{Z}/4$ of cyclic
permutations that naturally acts on $U$ and its subspaces $U_i$,
and turns $\Phi$ into a $\mathbb{Z}/4$-equivariant map. In turn
$\Phi_0$ is also a $\mathbb{Z}/4$-equivariant map and for the
proof of Theorem~\ref{thm:smooth-square-pegs} it would be
sufficient to show that such an equivariant map must have a zero.

Finally, let us record for the further reference that the
generator $\omega\in\mathbb{Z}/2$ acts on $S^1[4]=S^1\times W_4$
by reversing the orientation while the action on $U$ is the
antipodal action $\omega(v)=-v$, hence it preserves the
orientation of $U\cong \mathbb{R}^4$.

\subsection{Compactified configuration space $S^1[4]$}
\label{subsec:comp}

Let $S^1[4]$ be the canonical or \textit{FMASK}-compactification
of the configuration space $S^1(4)$. We use the basic properties
of this compactification, as outlined in Section~\ref{sec:fmask},
to define a modified test map ${\Psi}_0 : S^1[4]\rightarrow U$.

Let $\eta : (S^1)^4\rightarrow \mathbb{R}$ be the map defined on
the configuration $t=(t_1,t_2,t_3,t_4)\in (S^1)^4$ as the
arc-length diameter of the set $\{t_1,t_2,t_3,t_4\}$, i.e.\ the
minimum arc-length of a closed arc $L\subset S^1$ such that
$t_i\in L$ for each $i$. Let $\xi := \eta^{-1}$ and let $\Phi'$ be
the modification of the test map $\Phi$ (equation
(\ref{eqn:peg-test-map})) defined by
\begin{equation}\label{eqn:peg-test-map-modified}
\Phi'(t) := \xi(t)\cdot\Phi(t) =
(\xi(t)\phi_1(F(t)),\xi(t)\phi_2(F(t)),\xi(t)\phi_3(F(t))).
\end{equation}
Finally, let $\Phi_0'$ be the restriction of $\Phi_0$ on $S^1(4)$.

\begin{prop}\label{prop:extension}
The $\mathbb{Z}/4$-equivariant map $\Phi_0' : S^1(4)\rightarrow U$
can be extended to a $\mathbb{Z}/4$-equivariant map $\Psi :
S^1[4]\rightarrow U$ such that $\Psi(x)\neq 0$ for each $x\in
S^1[4]\setminus S^1(4)\cong S^1\times \partial W_4$. Moreover, the
$\mathbb{Z}/4$-equivariant homotopy class of the restriction
$\Psi^\partial : S^1\times\partial W_4\rightarrow U\setminus\{0\}$
does not depend on the embedding $f : S^1\rightarrow
\mathbb{R}^2$.
\end{prop}

\medskip\noindent
{\bf Proof:} The extension $\Psi$ is clearly unique (if it
exists). It is also clear that the only case to be discussed is
the case of points $q\in S^1[4]\setminus S^1(4)$ such that
$\eta(q)=0$, or equivalently $\xi(q)=+\infty$. These are the
points which corresponds to pentagons in
Figure~\ref{fig:cyclohedron} and can be characterized as limits in
$S^1[4]$ of sequences $q_n = \{t_i^n\prec t_j^n\prec t_k^n\prec
t_l^n\}$, where $(i,j,k,l)$ is a cyclic permutation of elements
$\{1,2,3,4\}$ and $\measuredangle(t_i^n,t_l^n)\mapsto 0$ as
$n\mapsto +\infty$. The last condition implies that all sequences
$(t_i^n)_{n=1}^{+\infty}$ converge to the same point $s\in S^1$.

The point $q\in S^1[4]$, which is the limit (in $S^1[4]$) of the
sequence $q_n\in S^1(4)$, is (in the language of
Section~\ref{sec:cyclohedron}) best visualized in the
$2$-dimensional screens described by equation
(\ref{eqn:2-screen}). Since the associated barycentric coordinates
$\frac{\measuredangle(t_i^n,t_j^n)}{\measuredangle(t_i^n,t_l^n)}$
etc.\ are all well defined as functions on $S^1[4]$, it remains to
be checked that the same applies to the functions
$\frac{\phi_i(F(t))}{\eta(t)}$ that appear in the test map
(\ref{eqn:peg-test-map-modified}),  i.e.\ that these quotient can
be meaningfully (and continuously) extended to points $q\in
S^1[4]$. Since in the small neighborhood of $s\in S^1$ the
function $f: S^1\rightarrow \mathbb{R}^2$ is approximated by a
linear function, i.e.\ the curve $\Gamma$ is in the vicinity of
$z=f(s)$ (up to a higher order infinitesimal) approximated by its
tangent line at $z\in \Gamma$, we make the following useful
observation.
\begin{enumerate}
\item[{$\mathbf O_1$}] The value of the test function
$\Psi=(\Psi_1,\Psi_2,\Psi_3)$  at a point $q\in S^1[4]$ is equal
to the value of the original test function $\Phi_0$ at an
``infinitesimal'' quadruple $q_n = \{t_i^n\prec t_j^n\prec
t_k^n\prec t_l^n\}$ approximating $q$, divided by the associated
``infinitesimal'' arc-length $\eta(q_n)$.
\end{enumerate}
It follows that $\Psi_1(q)$ is always a non-zero vector collinear
to the tangent vector of $\Gamma$ at $z=f(s)$ with the only
exception being the case of the point $q$ represented by an
``infinitesimal parallelogram'' i.e.\ if
\[
\lim_{n\mapsto\infty}
\frac{\measuredangle(t_i^n,t_j^n)}{\measuredangle(t_i^n,t_k^n)} =
\lim_{n\mapsto\infty}
\frac{\measuredangle(t_k^n,t_l^n)}{\measuredangle(t_j^n,t_l^n)} =
0.
\]
In this case it is not difficult to check that $\Psi_3(q)\neq 0$
which comletes the proof of the first part of the proposition.

For the second part, let us suppose that $f_0, f_1: S^1\rightarrow
\mathbb{R}^2$ are two smooth embeddings such that both maps
$s\mapsto df_i/ds, \, i=0,1$ have degree $+1$. Then the
independence of the $\mathbb{Z}/4$-homotopy class of the map
$\Psi^\partial$ from the embedding $f : S^1\rightarrow
\mathbb{R}^2$  follows from the fact that any two such embeddings
can be connected by a regular homotopy i.e.\ by a family $f_t,\,
t\in [0,1]$ of smooth embeddings such $df_t(s)/ds\neq 0$ for each
$s\in S^1$. \hfill $\square$

\subsection{The obstruction $\ldots$}
\label{subsec:pegs-obstruction}

It remains to be shown that no map in the
$\mathbb{Z}/4$-equivariant homotopy class of the map
$\Psi^\partial$ can be extended to $S^1\times W_4$, i.e.\ that
there does not exist a $\mathbb{Z}/4$-equivariant map ``?'' that
completes the square
 \begin{equation}\label{cd:obstruction}
\begin{CD}
S^1\times \partial(W_4)) @>\Psi^\partial>> S^3 \\
 @VVV @VV\cong V \\
S^1\times W_4 @>?>>
 S^3
 \end{CD}
\end{equation}
The obstruction to the extension problem (\ref{cd:obstruction})
lives in the equivariant cohomology group
$$
H^4_{\mathbb{Z}/4}((S^1\times W_4),S^1\times \partial(W_4));
\mathcal{Z}))
$$
where $\mathcal{Z} = \pi_3(U\setminus\{0\})\cong
H_3(U\setminus\{0\})\cong H_3(S^3)\cong \mathbb{Z}$ inherits the
$\mathbb{Z}/4$-module structure from the $\mathbb{Z}/4$-action on
$U$. By equivariant Poincar\' e duality this group is isomorphic
to the group
$$H_0^{\mathbb{Z}/4}((S^1\times W_4)\setminus S^1\times
\partial(W_4)); \mathcal{Z}\otimes \varepsilon))\cong \mathcal{Z}_{\mathbb{Z}/4}\cong
\mathbb{Z}/2$$ where $\varepsilon$ is associated orientation
character, i.e.\ the $\mathbb{Z}/4$-module $H_4(S^1[4],\partial
S^1[4]; \mathbb{Z})$.

\subsection{$\ldots$ and its evaluation} \label{subsec:trasverse}

We evaluate the obstruction in $\mathcal{Z}_{\mathbb{Z}/4}\cong
\mathbb{Z}/2$ by counting the zeros of a ``generic'' (transverse
to zero) map $? : S^1\times W_4\rightarrow U$ (diagram
(\ref{cd:obstruction})) which extends a map in the
$\mathbb{Z}/4$-equivariant homotopy class of $\Psi^\partial$.

\begin{figure}[hbt]
\centering
\includegraphics[scale=.50]{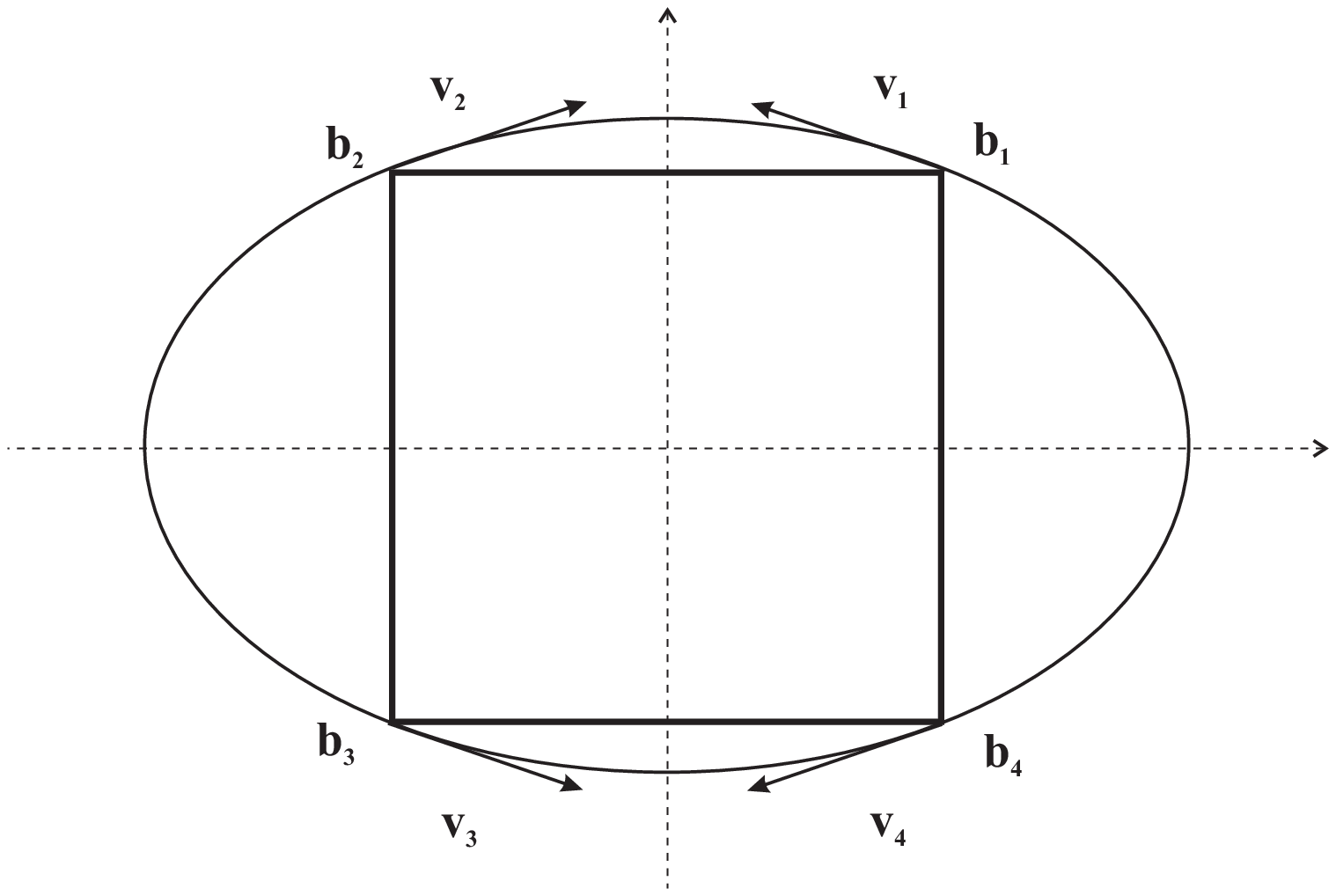}
\caption{} \label{fig:elipsa}
\end{figure}

Suppose that $\Gamma$ is a smooth oval in the plane which admits a
$\mathbb{Z}/2\times\mathbb{Z}/2$ symmetry. For example we can
choose for $\Gamma$ the ellipse centered at the origin, symmetric
with respect to the coordinate axes (Figure~\ref{fig:elipsa}).
Such an oval (ellipse) has a unique inscribed square. Suppose that
the vertices of this square (in counterclockwise order) are
$b_1,b_2,b_3,b_4$ and that $b_1$ is in the first quadrant.
Moreover we assume that $b_j = f(a_j)$ for some parameters $a_j\in
S^1$.

There is an obvious isomorphism $T_a(S^1(4))\cong
\oplus_{i=1}^4~T_{a_i}(S^1)$ of tangent spaces. Let $x_i$ be a
local coordinate on $S^1$ defined in the neighborhood of $a_i$.
For example let $x_i(c)$ be the (oriented) angle
$\measuredangle(a_i,c)$ swept by the radius vector moving from
$a_i$ to $c$. Let $[\frac{\partial}{\partial x_j}]_{j=1}^4
=[\frac{\partial}{\partial x_1}, \frac{\partial}{\partial x_2},
\frac{\partial}{\partial x_3}, \frac{\partial}{\partial x_4}]$ be
the associated basis (frame) of tangent vectors in $T_a(S^1(4))$.

We want to show that the differential $d\Psi_a :
T_a(S^1(4))\rightarrow T_0(\mathbb{R}^4)\cong \mathbb{R}^4$ of
$\Psi$, evaluated at $a=(a_1,a_2,a_3,a_4)$, is non-degenerate. Let
$y_i:=x_i\circ f^{-1}$ be the local coordinate on $\Gamma$ defined
in the neighborhood of $b_i$, induced by $x_i$. It follows that
the differential $dF_a : T_a((S^1)^4)\rightarrow T_b(\Gamma^4)$
maps the frame $[\frac{\partial}{\partial x_j}]_{j=1}^4$ to
$[\lambda_j\frac{\partial}{\partial y_j}]_{j=1}^4$, for
appropriate non-zero scalars $\lambda_j$.

Let $\alpha : (\mathbb{R}^2)^4\rightarrow
\mathbb{R}^2\times\mathbb{R}\times\mathbb{R}$ be the map defined
by $\alpha(y) = (\phi_1(y),\phi_2(y),\phi_3(y))$, so in particular
$\Phi(x) = \alpha(F(x))$.

The frame $[\frac{\partial \Psi}{\partial x_j}]_{j=1}^n =
[d\Psi(\frac{\partial }{\partial x_j})]_{j=1}^4$ is equal, up to
rescaling and possibly up to some changes of signs, to the frame
$[d\alpha (v_i)]_{i=1}^4$ where $v_i$ is an arbitrary (non-zero)
vector in $T_{b_i}(\Gamma)$ prescribed in advance. For convenience
(Figure~\ref{fig:elipsa}) we assume that the collection
$\{v_i\}_{i=1}^4$ is also $(\mathbb{Z}/2\times
\mathbb{Z}/2)$-invariant.

\begin{figure}[tbh]
\centering
\includegraphics[scale=.45]{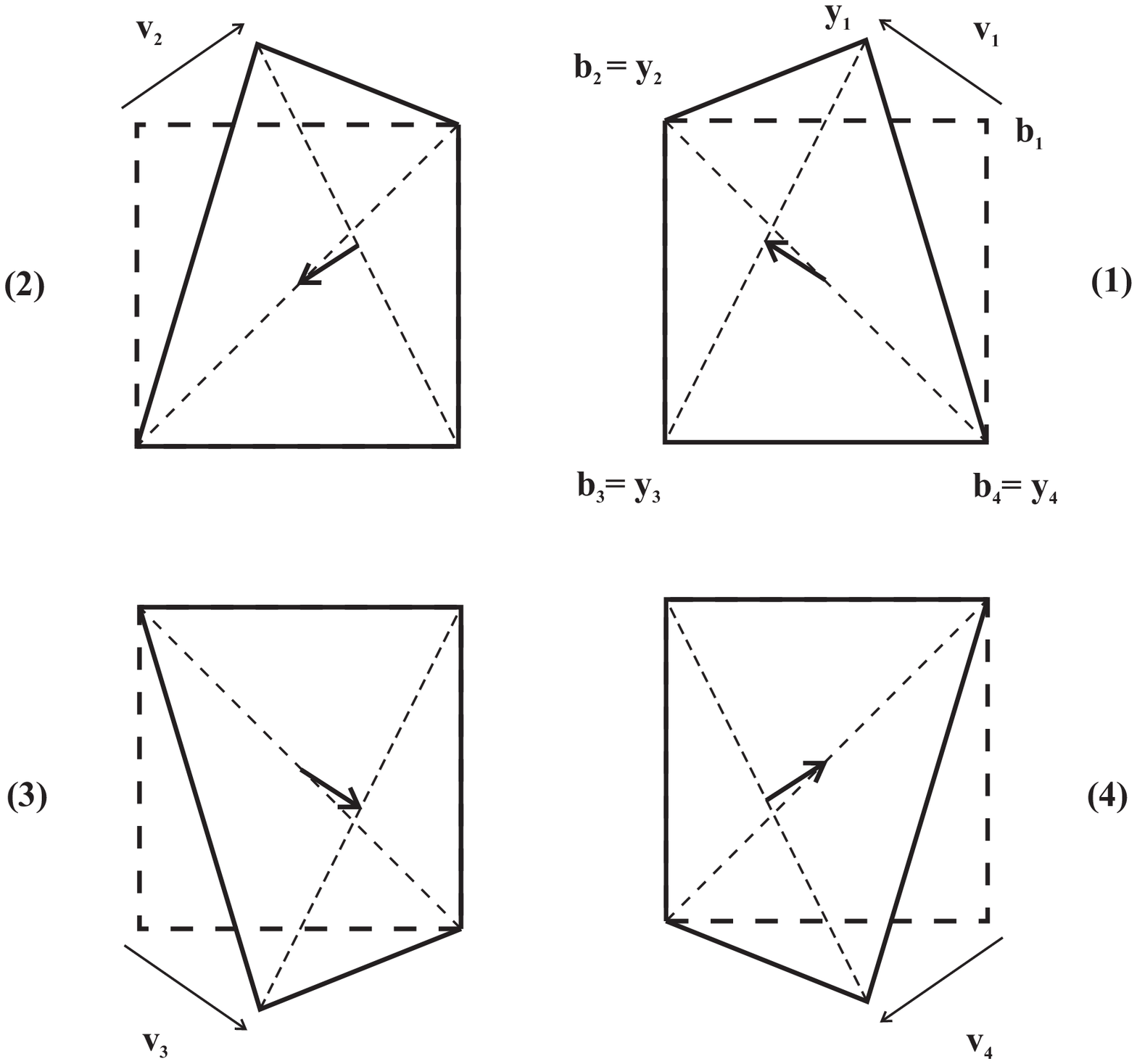}
\caption{} \label{fig:kvadrati}
\end{figure}

Let us suppose that the rate of change of $\alpha$ in the
direction of vector $v_1$, evaluated at the point
$(b_1,b_2,b_3,b_4)$,  is
 $$
d\alpha(v_1)=d\alpha_b(v_1)=(u;s,t)=(u_1,u_2;s,t)\in
\mathbb{R}^2\oplus \mathbb{R}\oplus \mathbb{R}.
 $$
By taking into account the $(\mathbb{Z}/2\times
\mathbb{Z}/2)$-symmetry of the curve $\Gamma$, one easily deduces
that
\[
\begin{array}{cccccccc}\label{eqn:derivatives}
d\alpha(v_2) &=& (u_1,-u_2;-s,t) &&& d\alpha(v_1) &=& (u_1,u_2;s,t)\\
d\alpha(v_3) &=& (-u_1,-u_2;s,t) &&& d\alpha(v_4) &=&
(-u_1,u_2;-s,t).
\end{array}
\]
It follows that the determinant $Det$ of the frame
$[d\alpha(v_j)]_{j=1}^4$ is
\[
Det = \left| \begin{array}{lrrr}
 u_1 & u_1 & -u_1 & -u_1 \\
 u_2 & -u_2 & -u_2 & u_2 \\
 s & \hspace{-5pt}-s & s & \! -s \\
 t & t & t & \! t
\end{array}
\right| = -16stu_1u_2\neq 0
\]
which in turn implies that the frame
$[\frac{\partial\Psi}{\partial x_j}]_{j=1}^4$ is also
non-degenerate.

\subsection{The proof of Theorem~\ref{thm:smooth-square-pegs}}
\label{subsec:proof-smooth-sqpegs}

\noindent {\bf Proof of Theorem~\ref{thm:smooth-square-pegs}:}
Assume that $f: S^1\rightarrow \mathbb{R}^2$ is a
(counterclockwise) smooth parametrization of the curve $\Gamma$.
The zeros of the associated ($\mathbb{Z}/4$-equivariant) ``test
map'' $\Phi_0: S^1(4)\rightarrow U$ (Section~\ref{subsec:config})
are in one-to-one correspondence with the squares inscribed in
$\Gamma$. After rescaling by a suitable positive, real function
$\xi$, the modified test map $\Phi_0':=\xi\Phi_0$ is
$\mathbb{Z}/4$-equivariantly extended (Section~\ref{subsec:comp})
to a map $\Psi : S^1[4]\rightarrow U$, where $S^1[4]$ is the
Fulton-MacPherson compactification of $S^1(4)$. The restriction
$\Psi^\partial$ of $\Psi$ on the boundary $\partial S^1[4]$ of
$S^1[4]$ has no zeros (Proposition~\ref{prop:extension}).
Moreover, its $\mathbb{Z}/4$-equivariant homotopy class is
independent of the original curve $\Gamma$. The obstruction for
extending $\mathbb{Z}/4$-equivariantly the map $\Psi^\partial :
\partial S^1[4]\rightarrow U\setminus\{0\}$ to $S^1[4]$ is found to be
non-trivial (Sections~\ref{subsec:pegs-obstruction} and
\ref{subsec:trasverse}) which finally implies that $\Phi_0$ must
have a zero in $S^1(4)$. \hfill $\square$

\section{Gr\" unbaum's conjecture}
\label{sec:Grunbaum}

B. Gr\" unbaum, see \cite{Gru} page 85, conjectured that every
Jordan curve in the plane contains the vertices of an
affine-regular hexagon. By definition a hexagon in the plane is
affine-regular if it is the image of a regular hexagon by an
affine automorphism of the plane.  In contrast to the square peg
problem, as emphasized in \cite{Gru}, the Gr\" unbaum's conjecture
has been opened even for the case of smooth curves. In this
section we establish this case of the conjecture and refer the
reader to Section~\ref{sec:conclusion} for a brief discussion how
the smoothness condition can be relaxed.

\begin{theo}\label{thm:gr} Every $C^1$-smooth, simple, closed curve in the plane
contains either the vertices of an affine-regular hexagon or six
collinear points which are the limit configuration of a convergent
sequence of affine-regular hexagons.
\end{theo}

\medskip The proof of Theorem~\ref{thm:gr} follows the
same scheme used in Section~\ref{sec:square-pegs} so we focus our
attention on differences and relevant calculations. By assumption
$\Gamma\subset \mathbb{R}^2$ is a simple, closed, $C^1$-smooth
curve in the plane, the last condition saying that a
parametrization can be chosen so that the curve has a non-zero,
continuously moving tangent vector.

In analogy with the square peg problem we choose for the {\em
configuration space} the Fulton-MacPherson compactification
$S^1[6]$ of the space $S^1(6)$ of all labelled, cyclically ordered
$6$-element subsets $\{t_1\prec t_2\prec\ldots\prec t_6\prec
t_1\}$ in $S^1$. Next we introduce the maps $\alpha, \beta,
\gamma, \delta$ which serve for testing if the points
$x_1,x_2,\ldots, x_6$ are consecutive vertices of an
affine-regular hexagon in the plane,

\begin{equation}\label{eqn:hexagon}
\begin{array}{l}
\alpha(x)= x_1 +x_4 - x_2 - x_5 \\
\beta(x)= x_2 + x_5 - x_3 - x_6 \\
\gamma(x) = x_3 + x_6 - x_1 - x_4\\
\delta(x) = x_1-x_2+x_3-x_4+x_5-x_6.
 \end{array}
\end{equation}
Note that the condition $\alpha(x)=0$ says that the midpoints of
the large diagonals $[x_1,x_4]$ and $[x_2,x_5]$ coincide while
$\delta(x)=0$, in addition to $\alpha(x)=\beta(x)=\gamma(x)=0$,
guarantees that the pairs of opposite sides are parallel to and
half the length of the large diagonal separating them.

\begin{figure}[htb]\vspace{-5cm}
\input{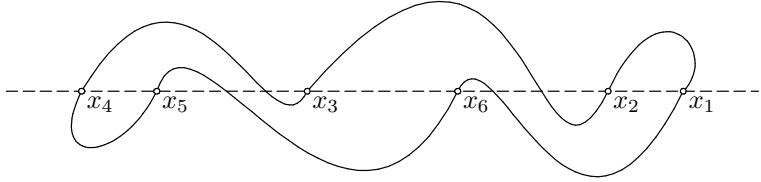}
\vspace{-2cm} \caption{A degenerate, affine-regular
hexagon.}\label{fig:deg-hexagon}
\end{figure}

As in the ``square peg problem'',  the system of equations
\begin{equation}\label{eqn:system}
\alpha(x)=\beta(x)=\gamma(x)=\delta(x)=0,
\end{equation}
admits, aside from genuine affine-regular hexagons, also some
degenerate solutions, for example the hexagons where all vertices
collapse to the same point. More generally there exist solutions
where points coincide in pairs, e.g.\ the solution $x_1=x_2=a,
x_4=x_5=-a, x_3=x_6$ and the solutions obtained from this one by a
cyclic permutation of indices. These two types of degenerate
solutions will be referred to as $1$-point and $3$-point
degenerate solutions. The system (\ref{eqn:system}) has also
collinear $6$-point solutions and together these are the only
degenerate solutions that can occur. In order to understand better
these $6$-point ``pseudo-solutions'', let us assume that
$x_1+x_3+x_5=x_2+x_4+x_6=0$ and that all these points belong to
the real axes. In light of the fact that $x_1+x_3=x_2,
x_2+x_4=x_3$, etc.\ we see that if $0<x_2<x_1$ then $x_3<0$
(Figure~\ref{fig:deg-hexagon}) which leads to the following simple
but important observation.

{\rm \begin{observ}\label{observ:hexagon} Suppose that
$x_1,x_2,\ldots, x_6$ are collinear points which are also vertices
of a degenerate, affine-regular hexagon, i.e.\ a configuration
obtained as a limit of a convergent sequence of affine-regular
hexagons. Suppose that these points appear in this order on a
smooth, Jordan curve $\Gamma$, i.e.\ $x_j=f(t_j)$ where $t_1\prec
t_2\prec \ldots \prec t_6 \prec t_1$. Then the order of the
appearance of these points on the line (in any direction) is a
non-trivial permutation of indices $1,2,\ldots, 6$ different from
a cyclic permutation.
\end{observ}}

\subsection{Compactified configuration space and the test maps}

 Let $U' \cong U_1^3\oplus U_2$, where $U_1\cong U_2\cong
\mathbb{R}^2$, be the preliminary test space defined as the total
target space for the test maps $\alpha,\beta,\gamma$ and $\delta$,
described in (\ref{eqn:hexagon}). Since the first three maps are
not independent $\alpha + \beta + \gamma = 0$, let
$V:=\{(u,v,w)\in U_1^3\mid u+v+w=0\}$ and let the actual ``test
space'' be the direct sum $U = V\oplus U_2\subset U'$. Let $F :
(S^1)^6\rightarrow (\mathbb{R}^2)^6$ be the map induced by the
embedding $f : S^1\rightarrow \mathbb{R}^2$ and $\Phi :
(S^1)^6\rightarrow U$ the ``test map'' where
$\Phi(t):=(\alpha(F(t)),\beta(F(t)),\gamma(F(t)))\oplus
(\delta(F(t)))$. Finally, let $\Phi_0 : S^1(6)\rightarrow U$ be
the restriction of $\Phi$ on the configuration space $S^1(6)$.

The next step, as in Section~\ref{subsec:comp}, is a rescaling of
the map $\Phi_0$ in order to make it suitable for an extension on
the compactified configuration space $S^1[6]$. As before
(Section~\ref{subsec:comp}) let $\eta : S^1(6)\rightarrow
\mathbb{R}^+$ be the map evaluating the ``circular diameter'' of a
configuration $t_1\prec t_2\prec\ldots\prec t_6\in S^1(6)$ and
$\xi:=\eta^{-1}$. Let $\Phi' := \xi\cdot \Phi$ be the rescaled
version of the map $\Phi$ and $\Phi_0'$ its restriction on
$S^1(6)$.

 The proof of the following proposition is similar
to the proof of Proposition~\ref{prop:extension} so we omit most
of the details.
\begin{prop}\label{prop:extension-2}
The $\mathbb{Z}/6$-equivariant map $\Phi_0' : S^1(6)\rightarrow U$
can be extended to a $\mathbb{Z}/6$-equivariant map $\Psi :
S^1[6]\rightarrow U$ such that $\Psi(x)\neq 0$ for each $x\in
S^1[6]\setminus S^1(6)\cong S^1\times \partial W_6$. Moreover, the
$\mathbb{Z}/6$-equivariant homotopy class of the restriction
$\Psi^\partial : S^1\times\partial W_6\rightarrow U\setminus\{0\}$
does not depend on the (counterclockwise) embedding $f :
S^1\rightarrow \mathbb{R}^2$.
\end{prop}

\medskip\noindent
{\bf Proof:} (outline) In order to show that $\Psi$ has no zeros
on the boundary $S^1\times \partial W_6$ we have to show that
``infinitesimal degenerate hexagons'' cannot appear as zeros of
the map $\Psi : S^1[6]\rightarrow U$. Indeed, this is ruled out by
the Observation~\ref{observ:hexagon}. The rest of
Proposition~\ref{prop:extension-2} is established by the arguments
already used in the proof of Proposition~\ref{prop:extension} so
we omit the details. \hfill $\square$

\subsection{The obstruction and its evaluation}
\label{subsec:obstruction-Grunbaum}

As in Section~\ref{subsec:pegs-obstruction}, there arises an
extension problem
\begin{equation}\label{cd:obstruction-Grunbaum}
\begin{CD}
S^1\times \partial(W_6)) @>\Psi^\partial>> S^5 \\
 @VVV @VV\cong V \\
S^1\times W_6 @>?>>
 S^5
 \end{CD}
\end{equation}
with the corresponding obstruction living in the group
$$
H^6_{\mathbb{Z}/6}((S^1\times W_6),S^1\times \partial(W_6));
\mathcal{Z}))\cong H_0^{\mathbb{Z}/6}((S^1\times W_6)\setminus
S^1\times
\partial(W_6)); \mathcal{Z}))\cong \mathcal{Z}_{\mathbb{Z}/6}\cong
\mathbb{Z}/2.
$$
As in Section~\ref{subsec:pegs-obstruction} we evaluate the
obstruction by choosing a conveniently ``generic'' curve
$\Gamma\subset \mathbb{R}^2$ and counting the number of
affine-regular hexagons inscribed in this curve.

\begin{figure}[htb]\vspace{-3cm}
\input{heksagon-1.pic}
\vspace{-3cm} \caption{}\label{fig:hexagon}
\end{figure}
Let $\Gamma$ be the boundary of a triangle
(Figure~\ref{fig:hexagon}). In order to make it a smooth curve one
is allowed to round its corners, however this will not affect the
calculations.

As it is clear from the picture there is only one affine-regular
hexagon inscribed in this curve. It remains to be shown, as in
Section~\ref{subsec:trasverse}, that a neighborhood of this
hexagon is mapped by the test map to a neighborhood of $0$, i.e.\
that $0$ is a regular point of the test map $\Psi$.

Assume that $x_1, x_2,\ldots, x_6$ are local coordinates (on the
configuration space $\Gamma(6)\cong S^1(6)$), in the neighborhood
of the hexagon depicted in Figure~\ref{fig:hexagon}. More
precisely $x_1$ is a point in the neighborhood of $x_1^0$,
constrained to move only on the $BC$ side of the triangle,
similarly $x_2,\ldots, x_6$ are perturbations of respective points
$x_2^0,\ldots, x_6^0$ allowed to move only on the boundary of the
triangle.

It follows, since the function $\gamma$ can be expressed in terms
of $\alpha$ and $\beta$, that we have to compute and establish the
non-triviality of the Jacobian $J$ of the map $x=(x_1,x_2,\ldots,
x_6)\mapsto (\alpha(x),\beta(x),\delta(x))$, evaluated at the
point $(x_1^0,\ldots, x_6^0)$.

By inspection, and up to some rescaling of vectors
$\overrightarrow{AB}, \overrightarrow{BC}, \overrightarrow{CA}$,
the Jacobian matrix is found to have the following form:
\begin{equation}\label{eqn:Jacobian-Grunbaum}
J =  \begin{array}{r|cccccc}
  & x_1 & x_2 & x_3 & x_4 & x_5 & x_6 \\ \hline \vspace{-3mm}&&&&&&\\
 \alpha & {\small\overrightarrow{BC}} & -{\small\overrightarrow{BC}} &
 0 & {\small\overrightarrow{CA}} & -{\small\overrightarrow{AB}} & 0 \vspace{1mm}\\
 \beta & 0 & {\small\overrightarrow{BC}} & -{\small\overrightarrow{CA}} &
 0 & {\small\overrightarrow{AB}} & {\small\overrightarrow{-AB}} \vspace{1mm}\\
\delta & {\small\overrightarrow{BC}} &
-{\small\overrightarrow{BC}} & {\small\overrightarrow{CA}} &
-{\small\overrightarrow{CA}} & {\small\overrightarrow{AB}} &
-{\small\overrightarrow{AB}}
\end{array}
\end{equation}
Finally, by choosing $\overrightarrow{BC}, \overrightarrow{CA}$
and $\overrightarrow{AB}$ to be respectively the column vectors of
the matrix
\[
\left[\begin{array}{crr} 1 & -1 & 0\\
0 & 1 & -1
\end{array}\right]
\]
we obtain a matrix with the determinant equal to $3$. This
calculation confirms that the hexagon depicted in
Figure~\ref{fig:hexagon} is indeed a non-degenerate solution of
the system of equations (\ref{eqn:system}) which in turn implies
that the obstruction to the extension problem
(\ref{cd:obstruction-Grunbaum}) is a non-trivial element of the
group $\mathcal{Z}_{\mathbb{Z}/6}\cong \mathbb{Z}/2$. \hfill
$\square$

\section{Hadwiger's conjecture}
\label{sec:Hadwiger}

\begin{conj}\label{conj:had}
{\bf (\cite{h71})} Every  simple closed curve in the Euclidean
$3$-space contains four distinct points which are the vertices of
a parallelogram.
\end{conj}

Relying on the same method as in the previous sections we
establish a stronger statement (at least for $C^1$-smooth curves)
that this parallelogram can be claimed to have all sides pairwise
equal i.e.\ to be a rhombus. As it turned out this theorem was
already established by Makeev in \cite{Mak} whose initial
motivation was the question of inscribing equilateral polygonal
lines in space curves.

\begin{theo}\label{thm:Hadwiger} {\rm (\cite{Mak})}
Every $C^1$-smooth simple closed curve $\Gamma$ immersed in the
Euclidean $3$-space contains the vertices of a rhombus.
\end{theo}

\medskip\noindent
{\bf Proof:} By assumption the curve $\Gamma\subset \mathbb{R}^3$
admits a $C^1$-parametrization. In other words it can be
parameterized by an injective $C^1$-mapping $\varphi : S^1 \to
\mathbb{R}^3$ such that the tangent vector $d\varphi/dt$ is
nowhere zero, continuous function of $t$.

As before $S^1(4)$ is the configuration space of cyclically
ordered four-tuples of distinct points on the circle. Given a
point $(t_1\prec t_2\prec t_3\prec t_4)\in S^1(4)$, let
$x_i=\varphi(t_i)$ for $i=1,2,3,4$. Let $\Gamma(4):={\rm
Image}(\Phi)$ where $\Phi : S^1(4)\rightarrow (\mathbb{R}^3)^4$ is
the map induced by $\varphi$.

Define a test map $F = F_\varphi : S^1(4)\to \mathbb{R}^4$ as the
composition $F:= \Psi\circ\Phi$ where $\Psi=(\Psi_1,\Psi_2):
(\mathbb{R}^3)^4\rightarrow \mathbb{R}^3\oplus \mathbb{R}$ is
described by equations:
\begin{equation}\label{eqn:test-map-Hawdiger}
\begin{array}{lcl}
\Psi_1(x)&=&x_1-x_2+x_3-x_4\\ \Psi_2(x)&=& \Vert x_1-x_2\Vert
-\Vert x_2-x_3\Vert +\Vert x_3-x_4\Vert -\Vert x_4-x_1\Vert.
\end{array}
\end{equation}

For a given input $x=(x_1,x_2,x_3,x_4)\in \mathbb{R}^4$, the first
function $\Psi_1(x)$ describes a point in $\mathbb{R}^3$ which is
equal to $0$ if and only if the mid-points of the diagonals of the
quadrilateral with the vertices $x_1,x_2,x_3,x_4$ coincide (i.e.\
if it is a parallelogram). If in addition the second function
$\Psi_2$ is equal to $0$, then such a parallelogram have all sides
equal (i.e.\ it is a rhombus). This shows that the rhombuses
inscribed in the curve $\Gamma$ correspond to zeros of the test
map $F=F_\varphi$.

Let $S^1[4]\cong S^1\times W_4$ be the Fulton-MacPherson
compactification of the configuration space $S^1(4)$. As in the
previous sections one can extend, after some rescaling, the
function $F$ to a function $\tilde F = \tilde{F}_\varphi:
S^1[4]\to \mathbb{R}^4$. More explicitly, if $\eta(t)$ is the
arc-length diameter of the set $t=\{t_1,t_2,t_3,t_4\}\subset S^1$
and $\xi(t):=\eta(t)^{-1}$ then $\tilde{F}$ is the extension of
the map $\xi\cdot F : S^1(4)\rightarrow \mathbb{R}^4$.

The group $\mathbb{Z}/4$ acts on the configuration space $S^1(4)$
by cyclic permutations and this action could be extended to its
compactification $S^1[4]$. Moreover, both the test map $F_\varphi$
and its extension $\tilde{F}_\varphi$ are
$\mathbb{Z}/4$-equivariant. The target space $\mathbb{R}^4\cong
U_1\oplus U_2$ naturally splits as the sum of a $3$-dimensional
and a $1$-dimensional, real $\mathbb{Z}/4$-representation.

\begin{prop}\label{prop:extension-3}
The restriction $\tilde{F}_\varphi^\partial$ of the map
$\tilde{F}_\varphi$ on the boundary $S^1\times \partial(W_4)$ of
$S^1[4]=S^1\times W_4$ has no zeros. Moreover, the
$\mathbb{Z}/4$-equivariant homotopy class of the restriction
$\tilde{F}_\varphi^\partial: S^1\times
\partial(W_4)\rightarrow \mathbb{R}^4\setminus\{0\}$ depends
neither on the curve $\Gamma$ nor on its $C^1$-parametrization
$\varphi$.
\end{prop}

\medskip\noindent
{\bf Proof:} Let $\tilde{F}_\varphi^1 : S^1\times
\partial(W_4)\rightarrow U_1$ and $\tilde{F}_\varphi^2 : S^1\times
\partial(W_4)\rightarrow U_2$ be the components of the map
$\tilde{F}_\varphi^\partial = (\tilde{F}_\varphi^1,
\tilde{F}_\varphi^2) : S^1\times
\partial(W_4)\rightarrow U_1\oplus U_2\cong \mathbb{R}^4$.
By definition (equation (\ref{eqn:test-map-Hawdiger}))
$\tilde{F}_{\varphi}^j,\, j=1,2$, is the extension of the map
$\xi\cdot(\Psi_j\circ\Phi)$.

For the majority of the points $q\in S^1\times \partial(W_4)$
already the function $\tilde{F}_\varphi^1$ is non-zero. More
precisely $\tilde{F}_\varphi^1(q)$ is zero only if $q$ is a
degenerate parallelogram i.e.\ if $q\in S^1\times (I_1\cup I_2)$
where $I_1=[(23)(41),(41)(23)]$ and $I_2=[(12)(34),(34)(12)]$ are
the corresponding edges of the cyclohedron $W_4$ depicted in
Figure~\ref{fig:cyclohedron}. If $q\in I_1\cup I_2$ then
$\tilde{F}_\varphi^2(q)\neq 0$, which completes the proof of the
first part of the proposition.

\medskip
For the proof of the second part of the proposition let us begin
with a simple observation that the $\mathbb{Z}/4$-equivariant
homotopy class of $\tilde{F}_\varphi^\partial$ does not depend on
the smooth reparameterization of the curve $\Gamma$. Indeed, such
a reparametrization $\alpha : S^1\rightarrow S^1$ defines a
nowhere zero vector field on $S^1$, which can be linearly
contracted to the zero vector field. Hence, there is a smooth
homotopy between $\varphi$ and $\varphi':=\varphi\circ\alpha$
which induces a $\mathbb{Z}/4$-equivariant homotopy between
$\tilde{F}_\varphi^\partial$ and $\tilde{F}_{\varphi'}^\partial$.
Similar argument can be used in the case when two curves (knots)
$\Gamma_\varphi$ and $\Gamma_\psi$ are in the same isotopy class,
i.e.\ if they represent the same knot.

\medskip
For the general case let us suppose that $\varphi$ and $\psi$ are
two $C^1$-smooth embeddings (knots) of $S^1$ in $\mathbb{R}^3$
which are not necessarily $C^1$-isotopic. A naive candidate for a
$\mathbb{Z}/4$-equivariant homotopy between
$\tilde{F}_\varphi^\partial$ and $\tilde{F}_{\psi}^\partial$ is
the linear homotopy
\begin{equation}\label{eqn:homotopy}
G : S^1[4]\times [0,1]\rightarrow \mathbb{R}^4
\end{equation}
defined by $G(q,t)=(G^1(q,t),G^2(q,t)):=
(1-t)\tilde{F}_\varphi^\partial + (1-t)\tilde{F}_\psi^\partial$.

\begin{enumerate}
\item[] \noindent {\bf Observation~1:} The second component $G^2$
of the linear homotopy is non-zero on $A=S^1\times (I_1\cup I_2)$.
Indeed, the signs of both $\tilde{F}_\varphi^1(q)$ and
$\tilde{F}_\psi^1(q)$ (for $q\in A$) are the same, since they
depend solely on the labelling of the vertices of the degenerate
parallelogram $q$.
\end{enumerate}
It follows from Observation~1 that it is sufficient to define a
(nowhere zero) homotopy $G^1(q,t)$ for $q\in D:= \partial
S^1[4]\setminus (S^1\times (I_1\cup I_2))$. If
$\bar{F}_\varphi^1,\bar{F}_\psi^1 : D\rightarrow S^2$ are the
normalized maps associated to $\tilde{F}_\varphi^1$ and
$\tilde{F}_\psi^1$, where by definition $\bar{F}_\varphi^1(q):=
\tilde{F}_\varphi^1(q)/\Vert\tilde{F}_\varphi^1(q) \Vert$ and
$\bar{F}_\psi^1(q):= \tilde{F}_\psi^1(q)/\Vert\tilde{F}_\psi^1(q)
\Vert$, then it is sufficient to define a
$\mathbb{Z}/4$-equivariant homotopy $\bar{G}^1: D\times
I\rightarrow S^2$ between these two maps.

\begin{enumerate}
\item[]\noindent {\bf Observation~2:} Both maps
$\bar{F}_\varphi^1,\bar{F}_\psi^1$ can be canonically and
$\mathbb{Z}/4$-equivariantly factored through the space
$\tilde{F}(S^1,2)$ defined in Example~\ref{exam:primer}
(Section~\ref{sec:fmask}). More precisely there exist a
``universal'' $(\mathbb{Z}/4\mapsto \mathbb{Z}/2)$-equivariant map
$\chi : D\rightarrow \tilde{F}(S^1,2)$ and
$\mathbb{Z}/2$-equivariant maps $\alpha, \beta :
\tilde{F}(S^1,2)\rightarrow S^2$ such that $\bar{F}_\varphi^1 =
\alpha\circ\chi$ and $\bar{F}_\psi^1=\beta\circ\chi$.
\end{enumerate}
\begin{figure}[hbt]
\centering
\includegraphics[scale=.50]{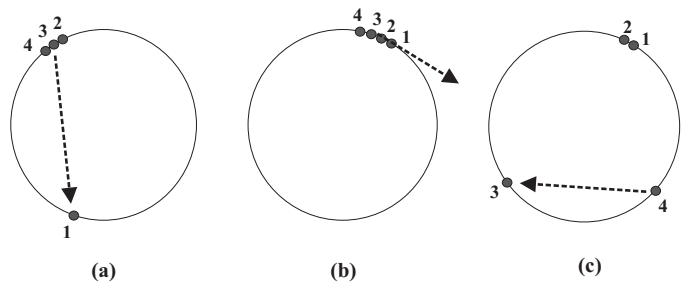}
\caption{The definition of the map $\chi$.}
\label{fig:3-circles-chi}
\end{figure}
The definition of the map $\chi$ is quite natural and motivated by
the definition of the map $\bar{F}_\varphi^1$ in the case when
$\varphi : S^1\hookrightarrow \mathbb{R}^2\subset \mathbb{R}^3$ is
essentially the identity map. Pictorially it is described in
Figure~\ref{fig:3-circles-chi}. Let $q_1=(x_1,x_2,x_3,x_4)$ be a
point in $\partial S^1[4]$. Generic examples are depicted in
Figure~\ref{fig:3-circles-chi} (where for simplicity $x_j$ is
labelled by $j$).

Recall that a point $\tilde{F}(S^1,2)$ is either an ordered pair
$(x,y)$ of points in $S^1$ or a pair $(z,u)$ where $z\in S^1$ and
$u$ is a unit tangent vector in $T_z(S^1)$. If
$q=(x_1,x_2,x_3,x_4)$ is the configuration depicted in
Figure~\ref{fig:3-circles-chi} (c), then $\chi(q):=(x_3,x_4)$. If
$q$ is a point depicted in Figure~\ref{fig:3-circles-chi} (a),
then $\chi(q):=(x_1,x_3)$. Finally if $q$ is a point depicted in
Figure~\ref{fig:3-circles-chi} (b), then $\chi(q)=(z,u)$ where
$z=x_1=x_2=x_3=x_4$ and $u$ is the unit tangent vector pointing in
he same direction as the vector $x_1-x_2+x_3-x_4$.

\medskip It is not difficult to show that $\tilde{F}(S^1,2)$ is,
as a $\mathbb{Z}/2$-space, $\mathbb{Z}/2$-homotopy equivalent to
the circle $S^1$ with the antipodal action.  The existence of the
$\mathbb{Z}/4$-homotopy $\bar{G}^1: D\times I\rightarrow S^2$
follows from Observation~2 and the fact that any two
$\mathbb{Z}/2$-maps $\alpha,\beta : \tilde{F}(S^1,2)\rightarrow
S^2$ are $\mathbb{Z}/2$-homotopic. \hfill $\square$

\begin{rem}\label{rem:Hadwiger}{\rm
The fact that the $\mathbb{Z}/4$-equivariant homotopy class of the
map $\tilde{F}_\varphi^\partial$ is the same, whether the curve
$\Gamma\subset \mathbb{R}^3$ is knotted or not, illustrates the
versatility of the ``cyclohedron approach'' to the problem of
finding polygonal pegs inscribed in space curves. This should be
compared to the planar case where any two (equally oriented)
embeddings are isotopic, hence the required homotopies can be
constructed by more direct methods. }
\end{rem}

\subsection{Obstruction and its evaluation}

{\bf Proof of Theorem~\ref{thm:Hadwiger} (cont.):}  Equipped with
Proposition~\ref{prop:extension-3}, we proceed as in the earlier
sections. The obstruction $O$ for a $\mathbb{Z}/4$-equivariant
extension of the (homotopically unique) map
$\tilde{F}^\partial_\varphi : S^1\times\partial W_4\rightarrow
\mathbb{R}^4\setminus\{0\}$ to $S^1\times W_4$ lives in
$H_{\mathbb{Z}/4}^4(S^1[4],\partial S^1[4];\mathcal{Z})$, cf.\
Section~\ref{subsec:pegs-obstruction}. The Poincar\' e dual of
this class belongs to the dual homology group
$H_0^{\mathbb{Z}/4}(S^1[4];\mathcal{Z}\otimes\epsilon)$ and can be
evaluated by a careful choice of the curve $\Gamma$.

\medskip
Let us consider the curve $\Gamma = \Gamma_1\cup J$ obtained as
the union of the non-closed curve $\Gamma_1={\rm Im}(\psi)$, where
$\psi : [0,2\pi]\to \mathbb{R}^3$ is defined by $\psi(t)=(\cos
t,\sin t, t)$, and the interval $J$ joining the endpoints
$(1,0,2\pi)$ and $(1,0,0)$ of the curve $\Gamma_1$. $\Gamma$ is
not a smooth curve however, by smoothing the corners, we obtain a
smooth curve $\Gamma'$ such that a rhombus $(x_1,x_2,x_3,x_3)$ is
inscribed in $\Gamma$ if and only if it is inscribed in $\Gamma'$.
For this reason we are allowed to use the test map associated to
the curve $\Gamma$.

Let $\varphi : [0,4\pi]\rightarrow \mathbb{R}^3$ be the
parametrization of $\Gamma$ defined by $\varphi(t) = \psi(t)$ for
$t\in [0,2\pi]$ and $\varphi(t):=(1,0,4\pi - t)$ for $t\in
[2\pi,4\pi]$.

We will show that there is a precisely one rhombus inscribed in
the curve $\Gamma$, i.e.\ a point $x=(x_1,x_2,x_3,x_4)\in
\Gamma(4)$ such that $F_\varphi(x)=0$. Moreover, it will be shown
that the test map $F=F_\varphi$ is transverse to $0\in
\mathbb{R}^4$, i.e.\ that $0$ is a regular value of $F$.

The only way to have two chords with coinciding mid-points are if
one chord has the end-points $(\cos t,\sin t,t)$ and $(\cos
(t+\pi),\sin (t+\pi),t+\pi)$ (and consequently the mid-point
$(0,0,t+\pi/2)$, where $0\leq t\leq \pi$), and the other chord has
the end-points $(\cos \pi,\sin \pi,\pi)=(-1,0,\pi)$ and
$(1,0,2t)$. Among the obtained parallelograms the only rhombus is
obtained when $t=\pi/2$.

Let us show that $0$ is indeed a regular value of the test map
$F=F_\varphi$. The vertices of the rhombus $q=(x_1,x_2,x_3,x_4)$
inscribed in the curve $\Gamma$ are:
\begin{equation}\label{eqn:vertices}
\begin{array}{lll}
x_1 = (0,1,\frac{\pi}{2}) &\quad & x_2 = (-1,0,\pi)\\
x_3 = (0,-1,\frac{3\pi}{2}) &\quad & x_4 = (1,0,\pi)
\end{array}
\end{equation}
The corresponding tangent vectors to the curve $\Gamma$ at these
points are:
\begin{equation}\label{eqn:tangent-vectors}
\begin{array}{lll}
\dot{x}_1 = (-1,0,1) &\quad& \dot{x}_2 = (0,-1,1)\\
\dot{x}_3 = (1,0,1) &\quad& \dot{x}_4 = (0,0,-1)
\end{array}
\end{equation}
By a slight abuse of language we can interpret
$\{\dot{x}_i\}_{i=1}^4$ also as a frame of tangent vectors at
$q=(x_1,x_2,x_3,x_4)\in \Gamma(4)$, i.e.\ as a basis of the
tangent space $T_q(\Gamma(4))$. Let us evaluate the rate of change
of functions $\Psi_1$ and $\Psi_2$ at $q\in \Gamma(4)\subset
\mathbb{R}^4$ in the directions of vectors $\dot{x}_i$, for
$i=1,2,3,4$. By definition (equation
(\ref{eqn:test-map-Hawdiger})) $\Psi_1(x)=x_1-x_2+x_3-x_4$, hence
\begin{equation}\label{eqn:rate-of-change}
\begin{array}{lll}
d\Psi_1(\dot{x}_1)=\dot{x}_1 = (-1,0,1) &\quad& d\Psi_1(\dot{x}_2)=-\dot{x}_2 = (0,1,-1)\\
d\Psi_1(\dot{x}_3)=\dot{x}_3 = (1,0,1) &\quad&
d\Psi_1(\dot{x}_4)=-\dot{x}_4 = (0,0,1).
\end{array}
\end{equation}
Let $\lambda=\sqrt{2+\pi^2/4}$ be the length of the side of the
rhombus $q=(x_1,x_2,x_3,x_4)$. Since by definition (equation
(\ref{eqn:test-map-Hawdiger}))
$$
\Psi_2(x)= \Vert x_1-x_2\Vert -\Vert x_2-x_3\Vert +\Vert
x_3-x_4\Vert -\Vert x_4-x_1\Vert
$$
we have
\begin{equation*}\label{eqn:rate-rate}
\begin{array}{l}
\lambda\, d\Psi_2(\dot{x}_1)=\langle x_1-x_2, \dot{x}_1 \rangle -
\langle x_1-x_4, \dot{x}_1 \rangle = \langle x_4-x_2, \dot{x}_1
\rangle = (2,0,0)\cdot (-1,0,1) = -2\\
\lambda\, d\Psi_2(\dot{x}_2)=\langle x_2-x_1, \dot{x}_2 \rangle -
\langle x_2-x_3, \dot{x}_2 \rangle = \langle x_3-x_1, \dot{x}_2
\rangle = (0,-2,\pi)\cdot (0,-1,1) = 2+\pi\\
\lambda\, d\Psi_2(\dot{x}_3)=\langle x_3-x_4, \dot{x}_3 \rangle -
\langle x_3-x_2, \dot{x}_3 \rangle = \langle x_2-x_4, \dot{x}_3
\rangle = (-2,0,0)\cdot (1,0,1) = -2\\
\lambda\, d\Psi_2(\dot{x}_4)=\langle x_4-x_3, \dot{x}_4 \rangle -
\langle x_4-x_1, \dot{x}_4 \rangle = \langle x_1-x_3, \dot{x}_4
\rangle = (0,2,-\pi)\cdot (0,0,-1) = \pi.
\end{array}
\end{equation*}
From here and equation (\ref{eqn:rate-of-change}) we conclude that
the Jacobian $dF$, evaluated at the point $t=\{t_1\prec t_2\prec
t_3\prec t_4\}$ parameterizing the rhombus $q$, is given by the
matrix
\[
\left[\begin{array}{rcrr} -1 & \, 0 & 1 & 0\\
0 & 1 & 0 & 0 \\ 1 & -1 & 1 & 1 \\
-2 & 2 + \pi & -2 & \pi
\end{array}\right]
\]
The determinant of this matrix is $-(2\pi + 4)\neq 0$ which
completes the proof that $0$ is a regular value of the test map
$F$.

The conclusion is that there is a precisely one rhombus inscribed
in the curve $\Gamma$ which is a regular value of the associated
test map. This implies that the obstruction element $O$ is
non-zero which completes the proof of Theorem~\ref{thm:Hadwiger}.
\hfill $\square$

\section{Concluding remarks and open problems}
\label{sec:conclusion}

\subsection{Relaxing the smoothness condition}\label{sec:extension}
The method of canonical compactifications was applied in
Sections~\ref{sec:outline}--\ref{sec:Hadwiger} to the problem of
inscribing the polygonal pegs in {\em smooth curves}. In this
section we briefly show how, with minimal modifications, the same
method can be extended and successfully applied to the case of
non-smooth curves satisfying some weaker, locally defined,
condition. As emphasized in Section~\ref{sec:outline}, the method
of canonical compactifications builds on the ``configuration
space/test map''-scheme, and introduces two important
modifications. The first modification applies to the configuration
space $S^1(n)$ which is extended (compactified) to the canonical
compactification (\textit{FMASK}-compactification) $S^1[n]$.
Secondly, the test map $\Phi : S^1(n)\rightarrow V$, for an
associated test space $V$, is modified to a new test map $\Psi:
S^1[n]\rightarrow V$. The local properties of the curve enter the
stage essentially in the following two ways.

\begin{enumerate}
\item[$\bullet_1$] The requirement that the Jordan curve is
$C^1$-smooth is essentially used
(Propositions~\ref{prop:extension}, \ref{prop:extension-2}, and
\ref{prop:extension-3}) in the definition of the modified test map
$\Psi: S^1[n]\rightarrow V$. \item[$\bullet_2$] Local properties
of the curve are used to guarantee that there are no
``infinitesimal squares'' inscribed in the curve which in turn
guarantees that the test map $\Psi$ has no zeros on the remainder
$S^1[n]\setminus S^1(n)$ of the compactification.
\end{enumerate}

The second condition is quite natural and in one form or another
it is present in all known results in this area. The most general
known conditions that rule out the existence of infinitesimal
inscribed squares  have been proposed by Stromquist \cite{St}. At
present it is not clear how this type of condition can be avoided
or at least considerably weakened.

Here we focus on the first condition $\bullet_1$ and show that in
principle one should be able to define the modified test map
$\Psi$ in all cases of interest, provided we are prepared to use
more general forms of \textit{canonical compactifications} which
include the \textit{FMASK}-compactification as a special case. In
other words one can always define the modified test map $\Psi :
\tau(S^1(n))\rightarrow V$, even if the curve $\Gamma$ is not
$C^1$-smooth, for a suitable compactification
$\tau(S^1(n))\geqslant S^1[n]$.

\begin{defin}\label{def:canonical}
Let $F[\mathbb{R}^2,n]$ be the canonical or
\textit{FMASK}-compactification of the configuration space
$F(\mathbb{R}^2,n)$ of all labelled, n-tuples of distinct points
in $\mathbb{R}^2$ (Definition~\ref{def:FMASK-2}). For a given
(oriented) Jordan curve $\Lambda\subset \mathbb{R}^2$ let
$\Lambda(n)$ be the collection of all $q=(q_1,\ldots, q_n)\in
F(\mathbb{R}^2,n)$ such that $q_i\in \Lambda$ for each $i$ and the
points $q_i$ appear on $\Lambda$ in the order which agrees with
the chosen orientation on $\Lambda$. The canonical
compactification $\Lambda[n]$ of $\Lambda(n)$ is defined as the
closure of $\Lambda(n)$ in $F[\mathbb{R}^2,n]$.
\end{defin}
The definition of $\Lambda[n]$ is in agreement with the definition
of $S^1[n]$ from Section~\ref{sec:fmask} provided $S^1$ is the
standard unit circle in $\mathbb{R}^2$. Canonical
compactifications $\Lambda[n]$ for a suitable $\Lambda$ can be
used as the source space for the test map $\Psi$. We illustrate
the key idea by giving some hints how the result about the square
pegs inscribed in Jordan curves can be established for {\em
piecewise smooth curves} $\Gamma$.

\medskip
Let $\Lambda_m$ be a regular $m$-gon in $\mathbb{R}^2$ (oriented
counterclockwise) and let $\Lambda_m[n]$ be the canonical
compactification of the configuration space $\Lambda_m(n)$. For
each piecewise smooth (oriented) Jordan curve $\Gamma\subset
\mathbb{R}^2$, which has at most $m$ non-smooth points, there
exists an (orientation preserving) homeomorphism $f :
\Lambda_m\rightarrow \Gamma$ which is smooth (with $df\neq 0$) on
each of the segments of $\Lambda_m$. Such a piecewise smooth
homeomorphism induces a continuous map $F :
\Lambda_m(n)\rightarrow \Gamma(n)$ which extends to the map
$\tilde{F}: \Lambda_m[n]\rightarrow \Gamma[n]$ of associated
canonical compactifications.

Moreover, the primary test map $\Phi : \Lambda_m(n)\rightarrow U$
defined by the equation (\ref{eqn:peg-test-map}) in
Section~\ref{sec:square-pegs}, on multiplication by the rescaling
factor $\eta$, can be extended to a test map $\Psi :
\Lambda_m[n]\rightarrow U$. This is established essentially by the
same argument already used in the proof of
Proposition~\ref{prop:extension}. The rest of the argument is
quite similar to the proof of the smooth case and does not require
new ideas.

\subsection{Problems and conjectures}\label{sec:open}

The following problem reflects our impression that the answer to
the Hadwiger's problem for smooth curves, given in
Section~\ref{sec:Hadwiger},  is somewhat exceptional. It seems
quite natural to expect that there should exist a space polygonal
peg of some sort which always appears in some knots while it is
missing in some realizations of other knots.

\begin{prob}\label{prob:knots}{\rm
Is there a genuine polygonal peg property that can distinguish
knots? In other words, is there a {\em naturally defined family}
$\mathcal{F}$ of polygonal pegs such that for some knot type
$\mathcal{N}_1$ knots $K\in \mathcal{N}_1$ always exhibit (grip) a
polygonal peg from the class $\mathcal{F}$ while for some other
knot type $\mathcal{N}_2$ there is a representative $L\in
\mathcal{N}_2$ which does not have this property.}
\end{prob}

If we stretch somewhat the concept of a ``naturally defined
family'' of polygonal pegs  by allowing families $\mathcal{F}$
that are purely non-metric in the sense that a polygonal peg $P$
belongs to $\mathcal{F}$ if and only if it satisfies a condition
based on (co)incidences of associated points and lines, then there
is a very interesting example showing that the answer to
Problem~\ref{prob:knots} should be affirmative. Indeed, as shown
in \cite{Pan} for generic polygonal knots and in \cite{Kup} for
tame knots (see also \cite{MM}), {\em quadrisecant lines} are
always present in nontrivial knots. On the other hand they
obviously may be absent in some realization of the unknot.
Moreover, it was shown in \cite{BSCS} that a weighted sum of
quadrisecants is a genuine knot invariant (the second coefficient
of the Conway polynomial). All this serves as a motivation for the
following bold conjecture.

\begin{conj}\label{conj:Vassiliev}
Polygonal pegs detect all finite type invariants.
\end{conj}

Conjecture~\ref{conj:Vassiliev} looks quite natural, however there
is an opposite point of view which deserves to be explored.
Suppose that the existence of a polygonal peg in a knot is
established by the \textit{CS/TM}-scheme, in the spirit of
Sections~\ref{sec:Grunbaum} and \ref{sec:Hadwiger}. The associated
test map incorporates a description of the polygonal peg in terms
of its characteristic metric and/or coincidence properties. For
example in the test map (\ref{eqn:test-map-Hawdiger}) for the
Hadwiger's problem (Section~\ref{sec:Hadwiger}) the first equation
tests the coincidence of mid-points of diagonals, while the second
equation records a pure metric property of a rhombus.

Suppose that a special test map for a concrete polygonal peg can
be designed so that it uses only the metric properties of the peg.
The reader is referred to \cite{Mak} and \cite{Mat}, Chapter~III,
for examples of such polygonal pegs.

Given a smooth knot $f : S^1\rightarrow \mathbb{R}^3$, one can
pull back the metric from the ambient space $\mathbb{R}^3$ to a
metric on $S^1$ and express the original polygonal peg problem as
a question about the existence of a peg in $S^1$, relative to this
metric. The punch line is that if a polygonal peg problem allows a
purely metric test map then, in light of the fact that any two
metrics on $S^1$ are homotopic, the associated obstructions should
be the same (cf.\ Remark~\ref{rem:Hadwiger}). As a consequence
such a polygonal peg problem cannot detect a knot.

\medskip\noindent
{\bf Acknowledgement:} It is a real pleasure to acknowledge
valuable comments and remarks by Benjamin Matschke, Igor Pak, and
Dev Sinha, as well as the hospitality of the Technical University
in Berlin (Discrete Geometry Group).


\begin{thebibliography}{abcd}

\bibitem[A-S]{AS} S.~Axelrod and I. Singer. Chern-Simons perturbation theory
II. \textit{Jour.\ Diff.\ Geom.} \textbf{39} (1994), no. 1,
173–-213.

\bibitem[B-T]{BT} R.~Bott and C.~Taubes. On the self-linking of
knots. \textit{J.\ Math.\ Phys.} \textbf{35} (1994), no.\ 10,
5247–-5287.

\bibitem[BSCS]{BSCS} R.~Budney, K.~Scannell, J.~Conant, and D.~Sinha.
New perspectives on self linking. \textit{Advances in
Mathematics}, Vol.\ 191 No 1 (2005), 78-113.

\bibitem[Die]{tomDieck} T. tom Dieck, \textit{Transformation Groups,}
De Gruyter studies in mathematics vol. 8, Berlin 1987.

\bibitem[Em]{Em} A.~Emch. Some properties of closed convex curves in the
plane. \textit{Amer.\ J.\ Math.}, \textbf{35} (1913), 407--412.


\bibitem[E]{Eng} R.~Engelking. \textit{General Topology}. PWN,
Warszawa 1977.

\bibitem[F-M]{FM} W.~Fulton and R.~MacPherson. Compactification of configuration
spaces. \textit{Ann.\ of Math.} \textbf{139} (1994), 183–-225.

\bibitem[Gri]{Gri} H.B.~Griffiths. The topology of square pegs in round
holes. \textit{Proc.\ London Math.\ Soc.} 62 (1991), 647–-672.

\bibitem[Gr\" u]{Gru} B.~Gr\" unbaum. \textit{Arrangements and
spreads}. AMS, Providence, RI, 1972.

\bibitem[Gug]{Gug} H. Guggenheimer. Finite sets on curves and
surfaces. \textit{Israel J.\ Math.} 3 (1965), 104–-112.


\bibitem[Ha]{h71} H.~Hadwiger, Ungel\"oste Probleme Nr. 53.
\textit{Elem.\ Math.} 26 (1971), 58.


\bibitem[HLM]{HLM} H.~Hadwiger, D.G.~Larman, and P.~Mani. Hyperrombs
inscribed to convex bodies. \textit{J.\ Combin.\ Theory Ser.\ B}
\textbf{24} (1978), 290--293.

\bibitem[Heb]{Heb} C.M.~Hebbert. The inscribed and circumscribed squares of a
quadrilateral and their significance in kinematic geometry.
\textit{ Ann.\ of Math.} 16 (1914/15), 38–-42.

\bibitem[Jer]{Jer} R.P.~Jerrard. Inscribed squares in plane curves.
\textit{Trans.\ AMS} 98 (1961), 234–-241.

\bibitem[Kak]{Kak} S.\ Kakeya. On the inscribed rectangles of a closed curvex
curve. \textit{T\^{o}hoku Math.\ J.} 9 (1916), 163–-166.

\bibitem[Kup]{Kup} G.~Kuperberg. Quadrisecants of knots and links. \textit{J.\
Knot Theory Ramifications}, 3 (1994) 41–50.

\bibitem[KlWa]{KlWa} V.L.Klee and S.~Wagon. \textit{Old and new unsolved
problems in plane geometry and number theory}. MAA, Washington,
DC, 1991.

\bibitem[Ko]{Ko} M.~Kontsevich. Operads and motives in deformation
quantization. \textit{Lett.\ Math.\ Phys.} \textbf{48} (1999),
35–-72.

\bibitem[Mar]{Ma} M.~Markl. Simplex, associahedron, and
cyclohedron. \textit{Higher Homotopy Structures in Topology and
Mathematical Physics}, Contemporary Math., vol.\ 227, Amer.\
Math.\ Soc., 1999, pp.\ 235--265.

\bibitem[MSS]{MSS} M.~Markl, S.~Shnider, and J.~Stasheff.
\textit{Operads in Algebra, Topology and Physics}, Math.\ Surveys
and Monographs 96, Amer.\ math.\ Soc., 2002.


\bibitem[Mak]{Mak} V.V.~Makeev. Quadrangles inscribed in a closed
curve and the vertices of a curve, \textit{J. Math. Sci. (N. Y.)},
Vol.~131, No.~1, 2005. Translated from \textit{Zap.\ Nauchn.\ Sem.
S.-Peterburg. Otdel. Mat. Inst. Steklov. (POMI))}.

\bibitem[Mat]{Mat} B.~Matschke. \textit{Equivariant Topologu and
Applications}, Diploma Thesis, TU Berlin, September 2008, \url{
http://www.math.tu-berlin.de/~matschke/DiplomaThesis.pdf}.


\bibitem[MM]{MM} H.R.~Morton, D.M.Q.~Mond. Closed curves with no
quadrisecants. \textit{Topology} 21 (1982) 235–243.


\bibitem[Pa08]{Pak} I.~Pak. The discrete square peg problem,
arXiv:0804.0657v1 [math.MG] 4 Apr 2008.

\bibitem[Pak]{Pa} I.~Pak. \textit{Lectures on Discrete and Polyhedral Geometry},
book in preparation, \url{http://www.math.umn.edu/~pak/book.htm}.

\bibitem[Pan]{Pan} E.~Pannwitz. Eine elementargeometrische Eigenschaft von
Verschlingungen und Knoten. \textit{Math.\ Ann.} 108 (1933)
629–672.


\bibitem[Si]{Si} D.~Sinha. Manifold-theoretic compactifications of con\-fi\-gu\-ra\-tion
spaces. \newline math.GT/0306385, 2003.

\bibitem[Shn]{Shn} L.G.~Shnirel'man. On some geometric properties
of closed curves (in Russian). \textit{Uspehi Matem.\ Nauk}
\textbf{10} (1944), 34--44; available at
\url{http://tinyurl.com/28gsy3}.


\bibitem[St]{St} W. Stromquist. Inscribed squares and square-like quadrilaterals in
closed curves. \textit{Mathematika} 36 (1989), 187–-197.

\bibitem[\v Z96]{User1}
R.~\v Zivaljevi\' c.
\newblock User's guide to equivariant methods
in combinatorics.  {\textit Pu\-bli\-cations de l'Institut
Mathematique} (Beograd), 59(73), 114--130, 1996.

\bibitem[\v Z98]{User2}
R.~\v Zivaljevi\' c.
\newblock User's guide to equivariant methods
in combinatorics II. {\textit Publi\-cations de l'Institut
Mathematique} (Beograd), 64(78) 1998, 107--132.

\bibitem[\v Z04]{Z04}
R.T. \v Zivaljevi\'{c}. Topological methods. Chapter 14 in
\textit{Handbook of Discrete and Computational Geometry}, J.E.\
Goodman, J.\ O'Rourke, eds, Chapman \& Hall/CRC 2004, 305 -- 330.




\end{thebibliography}
\end{document}